\documentclass{ametsoc}

\journal{jpo}

\bibpunct{(}{)}{;}{a}{}{,}

\title{Machine-Learning Ocean Dynamics from Lagrangian Drifter Trajectories}

\authors{Nikolas O. Aksamit \correspondingauthor{N.O. Aksamit, Swiss Federal Institute of Technology, Leonhardstrasse 21, 8092 Zurich, Switzerland}}

\email{naksamit@ethz.ch}

\affiliation{Institute for Mechanical Systems, Swiss Federal Institute of Technology, Zurich, Switzerland} 
\extraauthor{Themistoklis P. Sapsis}
\extraaffil{Department of Mechanical Engineering, Massachussetts Institute of Technology, Cambridge, Massachussetts, USA}
\extraauthor{George Haller}
\extraaffil{Institute for Mechanical Systems, Swiss Federal Institute of Technology, Zurich, Switzerland}

\abstract{Lagrangian ocean drifters provide highly accurate approximations of ocean surface currents but are sparsely located across the globe. As drifters passively follow ocean currents, there is minimal control on where they will be making measurements, providing limited temporal coverage for a given region. Complementary Eulerian velocity data is available with global coverage, but is itself limited by the spatial and temporal resolution possible with satellite altimetry measurements. Additionally, altimetry measurements approximate geostrophic components of ocean currents but neglect smaller sub-mesoscale motions and require smoothing and interpolation from raw satellite track measurements. In an effort to harness the rich dynamics available in ocean drifter datasets, we have trained a recurrent neural network on the time history of drifter motion to minimize the error in a reduced-order Maxey-Riley drifter model. This approach relies on a slow-manifold approximation to determine the most mathematically-relevant variables with which to train, subsequently improving the temporal and spatial resolution of the underlying velocity field. By adding this neural network component, we also correct drifter trajectories near sub-mesoscale features missed by deterministic models using only satellite and wind reanalysis data. The effect of varying similarity between training and testing trajectory datasets for the blended model was evaluated, as was the effect of seasonality in the Gulf of Mexico.}

\begin{document}

\maketitle

%

\section{Introduction}

Ocean currents play a critical role in the climate at timescales from weeks to decades. Ocean mixing is responsible for the scalar transport of temperature and salinity, impacting physical processes ranging from weather system formation, the melting rate of sea ice, the abundance and location of ocean ecosystems and the dispersion of microplastics \citep{Trenberth1994, Martin2002, Maximenko2012, Levy2018}. With the public availability of satellite altimetry measurements, researchers actively study global ocean currents and mesoscale features in the ocean, ranging in size from 10-200 km \citep{Stewart2008}, in near-real time from an Eulerian perspective. These data have allowed a better understanding of the role of prominent circulation features in the deep ocean, like the Gulf Stream and the Pacific Gyre, as well as smaller coherent structures like the Agulhas rings \citep{Wang2015}. In shallower waters, or at the interface of two mesoscale features, the influence of sub-mesoscale dynamics on ocean transport increases. The interaction of meso- and sub-mesoscale motions, however, is still poorly understood \citep{Beron-Vera2018}. Satellite observations typically represent the mixing and energy transport in sub-mesoscale features (1-10 km) inaccurately, with coastal features and tidal influences largely overlooked, as they do not have significant sea-surface-height signatures \citep{Thomas2008, Ferrari2008}.

Complementary Lagrangian measurements of ocean currents are also available, with their origins dating back to drifter deployments during the 1872 voyage of the HMS Challenger \citep{Niiler2001}. Since then, the quality and coverage of Lagrangian drifter data has improved through the use of satellite positioning systems and improved drifter designs. Lagrangian drifters come in varying geometries but are universally composed of a drogue that is fully submerged in the water, behaving as a sail, a float on the ocean surface that communicates location and velocity data, and a tether connecting the two. The size of the entire apparatus is typically on the order of 1-5 meters, with recent designs at the lower end of that range. The finite-size (inertial) effects and the buoyancy of every drifter prevent perfect adherence to fluid-particle pathlines. Considerable research has been conducted quantifying the ocean-drifter slip-velocity; thorough reviews of such processes have been presented over the years in great detail by \citet{Niiler1987}, \citet{Niiler1995}, \citet{DAsaro2003}, and others. Even with the disadvantages of non-zero slip velocities, Lagrangian drifters have been pivotal in our understanding of the oceans by providing velocity and position data for a range of motions spanning meters to thousands of kilometers. Such drifters continue to play a crucial role in understanding physical oceanography in changing climates with improving technology \citep{Lumpkin2016}. 

Typically, efforts to model Lagrangian drifter motions suffer in accuracy away from large structures, as noted in a recent evaluation of remotely-sensed products for drifter modeling near the Gulf Loop Current by \citet{Liu2014}. This increased error is, in large part, due to the relatively low temporal and spatial resolution available for ocean and wind velocity fields that drive drifter models. The focus of the current research is to harness the rich dynamics available in Lagrangian drifter data to enable better predictions of drifter velocities and trajectories from low-resolution Eulerian measurements of geostrophic ocean currents and wind reanalysis.

The present research capitalizes on theoretical work by \citet{Haller2008}, who developed reduced-order approximations of the computationally expensive Maxey-Riley equation for asymptotic inertial particle dynamics in unsteady flows. They identified an attracting slow-manifold in the equation and gave sufficient conditions for the manifold to attract other trajectories. This theoretical advance then allowed \citet{Beron-Vera2015} and \citet{Beron-Vera2016, Beron-Vera2019} to subsequently develop deterministic reduced-order models of ocean drifters that capture inertial aspects witnessed in cyclonic and anticyclonic eddies. In a numerical setting, \citet{Wan2018} and \citet{Wan2018a} also utilized slow-manifold-reduction to develop a leading-order particle motion model, but also introduced a blended model approach with a recurrent neural network to learn the higher-order terms in the slow dynamics. This blended machine learning and reduced-order model proved very successful in representing bubble motion for both laminar and turbulent flows while only using a limited laminar data set for neural network training \citep{Wan2018a}.

Here we use these recent modeling developments by first constructing a reduced-order model of Lagrangian ocean drifter motion similar to that of \citet{Beron-Vera2016}, but also including a-priori unknown functions for wind drag, wave rectification and near-surface shear forces. We then train a Long Short-Term Memory recurrent neural network \citep{Sepp1997} to learn the unknown functions and missing higher-order motion terms. Our slow-manifold approach provides us with an improved Lagrangian drifter model, and also determines the most relevant parameters and variables to use for learning sub-mesoscale dynamics. This development provides a new avenue for representing and understanding the role of sub-mesoscale motions in Lagrangian ocean dynamics and is fundamentally trainable for any basin worldwide. It is worth mentioning that the machine learning corrections developed likely depend strongly on the particular drifter geometry. Different drifters may not necessarily be characterized by the same data-driven model corrections just as different shapes will follow different trajectories in turbulent flow. As is shown below, the model is robust to different flow features, seasons and weather systems, and as with other machine learning applications, the similarity of trained and modeled features should be considered prior to widespread use. 


 \section{Methods}
 
Our main objective is to improve upon publicly-available geostrophic ocean current estimates for the Gulf of Mexico. Our approach fundamentally differs from data-assimilation methods by the inclusion of drifter history and special consideration of blended machine learning and deterministic models. We start by adding an Ekman current component induced by wind shear at the sea surface to geostrophic velocities obtained from satellite altimetry measurements. We then use this flow field in a reduced-order Maxey-Riley equation for baseline Lagrangian drifter-velocity predictions. A Long Short-Term Memory (LSTM) recurrent neural network \citep{Sepp1997} is then trained to minimize the modeled and real velocity residual along the real drifter trajectories. We evaluated the ability of our blended Maxey-Riley-LSTM model to predict drifter velocities given low-resolution input, as well as the ability to trace major flow features in a chaotic flow regime. Details of each of these steps can be found in the following subsections. 

\subsection{Data}

\subsubsection{Fluid Flow Fields}

The geostrophic velocity component, $\mathbf{u_g}$, used in our ocean current analysis is the freely available Archiving, Validation and Interpretation of Satellite Ocean Data (AVISO) L4 gridded velocity field hosted by the Copernicus Marine Environment Monitoring Service. The geostrophic current is the result of balancing pressure, as measured through sea surface height by satellite altimetry, and the Coriolis effect of a spinning earth. The equations of motion of the geostrophic flow are:

\begin{equation}
\frac{1}{\rho} \frac{\partial p}{\partial x} = f u_{g,2}, \qquad \frac{1}{\rho} \frac{\partial p}{\partial y} = -f u_{g,1},  \qquad  \frac{1}{\rho} \frac{\partial p}{\partial z} = -g, 
\end{equation}

where $\rho$ is the density of water, $p$ is the pressure, $g$ is the constant of gravity, $f$ is the Coriolis parameter and $\mathbf{u_g}=(u_{g,1}, u_{g,2})$.

We modify this flow field with an Ekman wind correction derived from ERA-5 reanalysis of 10-meter wind fields \citep{C3S}. Several recent ocean drifter modeling studies have shown that the addition of windage influence in ocean current estimation can greatly improve modeled drifter velocities \citep{Liu2014, Beron-Vera2016, Beron-Vera2019}. The ERA-5 data is freely available from the European Centre for Medium range Weather Forecasting (ECMWF) and is currently hosted by the Copernicus Programme Climate Change Service. We estimate the wind-shear contribution to ocean current using the Ekman model of wind influence, such that our resultant fluid velocity is $\mathbf{u}=\mathbf{u_g} + \mathbf{u_e}$, where

\begin{equation}
 \mathbf{u_e} = \left( u_{10}\frac{0.0127}{\sqrt{\text{sin}(\phi)}}  \text{cos}\left(\theta - \frac{\pi}{4}\right), \quad  u_{10}\frac{0.0127}{\sqrt{\text{sin}(\phi)}} \text{sin}\left(\theta - \frac{\pi}{4}\right) \right)
\end{equation}

$u_{10}$ is the 10 m wind speed, $\theta$ is the drifter direction in radians, and $\phi$ is the latitude \citep{Stewart2008}. More detailed estimations of wind drag are possible at this stage in the model development, but we instead rely on machine learning to account for slight inaccuracies in this estimation instead of further parameterization. This approach to including wind stress influences is comparable to the Geostrophic and Ekman Current Observatory (GEKCO) velocity products \citep{Sudre2013}. For our blended model, we used a spline interpolant to upscale the geostrophic velocity from the original 24-hour timestep and $1/4$ degree resolution to the 15-minute resolution of the drifter experiments. Additionally, the ERA-5 reanalysis was originally collected at a 6-hour timestep and $1/4$ degree resolution, and interpolated to the drifter resolution.

\subsubsection{Ocean Drifter Data}

The Lagrangian ocean drifter data used for training and validation in the current study were obtained from deployments during the LAgrangian Submesoscale ExpeRiment (LASER) \citep{Novelli2017, DAsaro2018, Haza2018} in January-February 2016 and the Grand LAgrangian Deployment (GLAD) in July 2012 \citep{Olascoaga2013, Poje2014, Jacobs2014, Coelho2015}. Position and velocity data during these extensive field campaigns in the Gulf of Mexico were measured onboard biodegradable CARTHE drifters \citep{Novelli2017} and were interpolated to 15-minute intervals. The GLAD experiment was conducted near the Deepwater Horizon spill site in the Gulf of Mexico during summer, with the follow-up LASER experiment, designed to identify sub-mesoscale motions in the same and different regions of the Gulf during winter. These experiments provide an unprecedented density of Lagrangian drifter trajectories near very complex flow features that are not clearly resolved in satellite altimetry or wind reanalysis data. The LASER and GLAD experiment data are a perfect test for assessing the feasibility of machine learning to adapt a reduced-order, low-resolution drifter model to complex flow features that are only resolved at high spatial and temporal resolution.

During LASER and GLAD, more than 1,100 drifters were deployed in the Gulf of Mexico, with the GLAD deployment being replicated by the first LASER deployment, and subsequent LASER deployments located near various structures in the Gulf. Many drifters were released on the continental shelf, close to the outlet of the Mississippi river in a region of strong salinity and temperature gradients at the freshwater-saltwater interface. Other drifters were deployed in deeper waters, thus providing drifter trajectories that measure sub-mesoscale motions of varying physical origin. Throughout the LASER experiment, several large El Niño winter storms occurred with significant surface wave generation, resulting in considerable drifter motion and sporadic drogue loss \citep{Haza2018}. As this loss is noted in the published dataset, we focused our training and modeling on the drogued drifters. The summertime GLAD deployment was affected by much stronger onshore-offshore winds than LASER and the velocity data contained much stronger evidence of windage and inertial oscillations \citep{Beron-Vera2015}.

\subsection{Blended Reduced-Order Data-Driven Drifter Model}

Starting from the Maxey-Riley (MR) equation of motion for a small, perfectly spherical particle in an unsteady flow \citep{Maxey1983}, we present a modified version for the motion of our Lagrangian drifter in an ocean current in the form:

\begin{equation}
\begin{split}
\rho_p \dot{\mathbf{v}} = &\rho_f \frac{D\mathbf{u}}{Dt} \\
& +(\rho_p -\rho_f) {\bf g} \\
&-\frac{9 \nu_f \rho_f}{2a^2} (\mathbf{v} - \mathbf{u} - \frac{a^2}{6}\Delta \mathbf{u}) \\
&- \frac{\rho_f}{2}(\dot{\mathbf{v}} - \frac{D}{Dt}(\mathbf{u}+\frac{a^2}{10}\Delta\mathbf{u})) \\
& -\frac{9 \nu_f \rho_f}{2a} \sqrt{\frac{\nu_f}{\pi}}\int_0^t\frac{1}{\sqrt{t-s}}(\dot{\mathbf{v}}(s)-\frac{d}{ds}(\mathbf{u}+\frac{a^2}{6}\delta\mathbf{u}))ds \\
& - \nu_a\rho_a\alpha(\mathbf{v} - {\bf u_{wind}}) - f(\rho_p \mathbf{v} - \rho_f \mathbf{u})^{\perp} +F({\bf u_ e},\mathbf{u}) .
\end{split}
\end{equation}

The first five terms on the right-hand side of the second equation in (3) represent the force of the current on the drifter, the buoyancy effect, the Stokes drag on the drogue, the force from the fluid moving with the particle, and the Basset-Bousinesq memory term, respectively. Specific to the motion of ocean drifters, we identified several additional forces that contribute significantly to a slip-velocity between the drifter and ocean current. These last forces include the Stokes drag of the surface wind on the exposed drifter float, the Coriolis force, the drag on the tether and drogue induced by shear currents, and motions induced by surface waves, \citep{Geyer1989, Niiler1995, Beron-Vera2015}. After combining the shear current and surface wave effects into one unknown function of current velocity and wind speed, $F({\bf u_ e},\mathbf{u})$, we are left with the remaining three terms in Eq 3. Following the slow-manifold approach described in the Appendix, we reduce this set of equations to the leading-order approximation

\begin{equation}
\dot{\mathbf{x}} = \mathbf{u}(\mathbf{x},t)+\epsilon\Big [ \left(1-\frac{3R}{2}\right)\left(\frac{D\mathbf{u}(\mathbf{x},t)}{Dt} - {\bf g}\right) -Rf(\delta_p-1)\mathbf{u}^{\perp}+H({\bf u_e}(\mathbf{x},t),\mathbf{u})\Big ] ,
\end{equation}

where $H$ is assumed to be a bounded function of the Ekman velocity and the current velocity, representing wind-wave and shear drag contributions; there is no known a-priori differential operation on $\mathbf{u}$, or $\mathbf{u}_e$ in $H$. For example, similar slip velocity estimates have been found to scale as linear functions of the ocean surface velocity \citep{Geyer1989, Edwards2006}. We avoid the assumptions necessary for such a linear slip-velocity estimate and relegate the learning of potentially more complicated wind-influence and shear drag terms to the data-driven component of our model. As detailed in the Appendix, the missing higher-order terms are shown to be functions of $\mathbf{u}$, $\mathbf{u}_e$, and $D\mathbf{u}/Dt$.

From this reduced-order model, we follow an approach similar to that of \citet{Wan2018a} by employing an LSTM to learn the contribution from the missing higher-order terms and the wind-wave function $H$. During our model order-reduction in the Appendix, we identify which variables the higher-order terms are functions of. This enables us to infer the time series with which to meaningfully train our recurrent neural network. We write out the 2D blended drifter model in the following form, neglecting vertical motion:

\begin{equation}
\begin{split}
\dot{{\bf x}} &=  {\bf v}({\bf x},t) + \epsilon\left[\left(\frac{3R}{2} - 1 \right)\frac{D{\bf u}}{Dt} -Rf(\delta_p-1){\bf u}^{\perp}\right]+\mathbf{G}(\boldsymbol{\xi}(t),\boldsymbol{\xi}(t-\tau),\boldsymbol{\xi}(t-2\tau), \dots) \\
 \boldsymbol{\xi}(t) &=  \left[{\bf u}({\bf x},t),{\bf u_e}({\bf x},t),  \frac{D{\bf u}}{Dt} \right]
 \end{split}
\end{equation}

Here $\mathbf{G}$ is a LSTM neural network, and $\boldsymbol{\xi}$ is a vector of time-delayed values of $\mathbf{u}$, $\mathbf{u}_e$ and of the material derivative $D\mathbf{u}/Dt$. Using our training set of real drifter trajectories, $\mathbf{G}$ is designed to learn both the error in our drifter approximation in Eq. 4 and the unknown function $H$, as a function of $\boldsymbol{\xi}$ over the full drifter history.

\section{Results}

We tested the blended model approach with a variety of training data sets from the LASER and GLAD experiments. The two experiments were similar, differing mainly in season (winter \& summer), the number of drifters and some of the seeding locations. Both experiments used the same CARTHE drifter design \citep{Novelli2017}. We trained and tested blended models on each season separately, as well as training a blended model with a combined set of drifters from both experiments. As mentioned previously, the combined drifter set includes trajectories influenced by varying weather systems with different magnitudes of inertial oscillation, thus displaying much more complex physics for the neural network to interpret.

Many of the drifters in both experiments were deployed in small clusters. By randomly selecting drifters, we created a rich training set that included many trajectories similar to those found in the test set. As machine learning approaches are most accurate when interpolating between features that the networks have already seen, using randomized training-test sets provided an easier task when testing the blended model. A cartoon example of random training and test drifter choice is shown in Figure 1a. However, with the currently limited amount of high-resolution data available for wider applications, we also tested the extrapolative ability of the blended model by limiting the similarity between training and test sets. We achieved this by handpicking the drifters, making sure drifters from the same tight deployment clusters were not separated between the training and test sets. This approach is demonstrated in contrast in Figure 1b. Because a minimal number of drifters is required for training, as well as a sufficient number of drifters for testing, this separation was imperfect. We will refer to these two methods of separating training and test drifters as random and non-repetitive from hereon. A comparison of the modeling capabilities of each method is discussed in the context of statistical measures of similarity in Sections 3.2 and 4.

The LASER campaign involved a much larger deployment of drifters than GLAD and hence serves as the primary source of discussion for the specific abilities of the blended model. We also briefly present the results from the summer GLAD model to complement the conclusions from the LASER data. By training and testing on a combined LASER and GLAD dataset, we show that with additional considerations, it may also be possible to train blended models for a progressively wider range of meso and sub-mesoscale features.

\subsection{Interpolative LASER Testing}

To test the ability of our blended model to recognize and predict motion around similar sub-mesoscale features, we initially trained the data-driven component $\mathbf{G}$ in Eq. 5 with 668 randomly selected Lagrangian drifter trajectories from the winter (LASER) dataset with time series length ranging from 16 to 60 days. The remaining 167 drifters were reserved as a test dataset in order to evaluate the model on new data. The results presented here come from using a simple LSTM architecture similar to that used by \citet{Wan2018a}. This architecture consists of one LSTM layer with 200 hidden units, and one fully connected layer trained for 100 epochs, and a time delay of $\tau$=15 minutes. This choice of training parameters has not been fully optimized, though we found it to result in sufficiently accurate predictions without overfitting when compared with training between 50 and 400 epochs and using LSTM layers of 50 to 400 units. For model comparison, we used our Ekman-modified flow field and the MR drifter model of \citet{Beron-Vera2015} as a baseline. The Beron-Vera et al. model is precisely Eq. 4 without the wind-wave function $H$.

After training, we evaluated the blended model’s single-step prediction ability of drifter velocity, that is, the model’s ability to predict velocity at each timestep along the actual drifter trajectory. As shown in Figure 2, the randomly-trained blended model showed considerable improvement in accuracy over the baseline model. Figure 2a displays histograms of the root mean squared error (RMSE) of the zonal velocity for the MR model and the blended model along all trajectories in the randomized LASER test set. A significant shift towards smaller errors with the blended model indicates an improvement in the model’s ability to accurately predict the 15-minute drifter velocity from daily $1/4$ degree geostrophic velocity and 6-hour wind reanalysis data. Example time series of the zonal component of modeled and actual drifter velocities are shown in Fig. 2b-d, with the ocean current  ${\bf u}({\bf x},t)$ in blue, output of the blended model in black and the real drifter velocity in orange. The distance to the nearest drifter in the training data set, the time to the soonest deployment of a trained drifter, and the RMSE of the blended model are all noted on the respective subplots.

We also tested the model’s multi-step predictive capabilities to evaluate how accurately we can forecast the trajectory of a real drifter in an environment with rich sub-mesoscale features. For multi-step prediction, we used blended model velocity output to update drifter positions over time with a fourth-order Runge-Kutta scheme, allowing modeling errors to accumulate in both the neural network memory and in the drifter position. For a chaotic flow such as in the ocean, it is nearly impossible for a perfect one-to-one model-drifter trajectory match. To account for uncertainty in the underlying chaotic dynamics, we introduced uncertainty in the initial conditions by seeding synthetic drifters in a small neighborhood ($<$10 km) of a real drifter position. This ensemble approach provided an insightful evaluation as dominant flow features began to appear in the cluster trajectories.

One test set example of this ensemble trajectory model is shown in Figure 3. The real LASER drifter trajectory is mapped as a bold red line in both the left and right subplots. The left panel shows the trajectories of all blended model drifters initialized in a small neighborhood in black, with their median trajectory displayed in yellow. The right panel shows analogous trajectories with blue Maxey Riley model output and the green median path. Over nearly 50 days of transport, the blended model drifters were able to trace out the large spiraling flow feature that was largely missed by the purely-deterministic MR equations.

We quantified the ability of the ensemble models to approximate real trajectories with the normalized cumulative separation distance, or skill-score, s from \citet{Liu2011}. For a drifter time series of length $n$, this metric normalizes the cumulative separation distance between the modeled drifter position $\mathbf{y}(t_j)$ and the actual drifter position $\mathbf{x}_D(t_j)$ by the cumulative length of the real drifter trajectory, $l(t_j)$, in the following form:

\[
c = \frac{\sum_{j=1}^{n}\Vert \mathbf{y}(t_j) - \mathbf{x}_D(t_j) \Vert}{\sum_{j=1}^{n}l(t_j)}, \qquad s=1-c.
\]

A skill score $s<0$ indicates the modeled and real drifter paths are deviating faster than the real drifter is actually moving. As shown in Fig. 4a, the mean skill score for an ensemble of blended drifters was better than the reduced-order model 67\% of the time. This improved trajectory modeling is also apparent in the skill score of the median ensemble trajectory (yellow and green paths in Fig. 3) in Fig. 4b, where the blended drifters outperformed the reduced-order model 68\% of the time. Of particular note are the different ranges of skill score values for the MR and blended model outputs. Specifically, there are large negative values present in the MR model evaluation that are non-existent for the blended model.

We investigated these negative skill scores and found that 75\% of them correspond with drifters deployed on Feb 7, 2016. In Figure 5a-b , we show the Feb 7 location of drifters near this deployment, with drifters colored by their skill scores for the blended and MR models, respectively. These colored coordinates are overlaid on the daily sea surface temperature from Aqua MODIS thermal infrared measurements on Feb 7. The MR drifters with the worst MR skill scores are those intentionally clustered in the strong mixing regions \citep{DAsaro2018} between the cold coastal and warm offshore waters in Fig. 5b. From this position, the warm waters extend south and east throughout the rest of the gulf and into the gulf loop current. For the blended model, there was no discernible difference in model performance across this interface, or elsewhere in the domain.

In Figure 5c-d, we show the median path drawn by 37 blended and MR drifters residing on the freshwater-saltwater interface in Fig. 5a-b. The MR ensemble exits the cool waters and quickly becomes entrained by the mesoscale geostrophic motions away from the coast. In contrast, the real drifters and blended model drifters stay near the mixing interface, in the cooler waters in the northern part of the domain. This sub-mesoscale model improvement was precisely the result of training the blended model on drifter velocity and acceleration data from similar mixing regions. The LSTM component was able to identify a signature of this mixing in $\boldsymbol{\xi}$ and prevent entrainment into the larger currents, forcing drifters to remain in the high temperature gradient mixing region on the shelf.

RMSE histograms of the random training of GLAD and LASER+GLAD datasets are shown in the following section. Skill scores from GLAD and LASER+GLAD datasets can be found in the supplemental information. Randomly selecting training drifters always resulted in model improvements over the MR equations. However, we analyze the robustness of these data-driven improvements in the following section by looking at minimally repetitive training and test sets.

\subsection{Extrapolative Drifter Testing}

We compared the blended model’s ability to correct velocities and trajectories for non-repetitive datasets from the LASER experiment, GLAD experiment, and a combined GLAD+LASER dataset to the randomized training method. In hopes to understand the different model performance with each training set, we generated several metrics comparing trajectories in test and training sets. To quantify how similar the training and test datasets were in the time domain, we calculated the maximum normalized correlation for zonal velocity of each test drifter against all drifters in the training set. To complement these values in the frequency domain, we also calculated the maximum mean magnitude spectral coherence between each test drifter and the training datasets. Lastly, we calculated distance of the initial position for each drifter to the nearest training drifter, and the time to the closest training drifter. Mean values of all these metrics from the randomized data separation are displayed in Table 1 for the LASER and GLAD data, with values inside parentheses indicating the non-repetitive separation scheme.

As can be seen in Table 1, maximum correlation and magnitude squared coherence change very little for the LASER data between the random and non-repetitive training schemes, but there is a significant increase in both the mean distance and time to nearest trained drifter. For the GLAD dataset, there was a decrease in correlation and MS coherence as well as an increase in distance and time to trained drifters when using the non-repetitive separation. To highlight the details of this influence, the RMSE values for single-step predictions of zonal velocities from the random and non-repetitive LASER models are compared to correlation and MS coherence in Fig. 6. There appears to be little relationship between similar test and training power spectra and low RMSE values, though a slight decrease in RMSE can be noted with increasing correlation, especially for the random training set. The red dot indicates a drifter present in both the random and non-repetitive test sets that is further analyzed in section 4.

The overall influence of the decreased similarity between trained and tested trajectories can be found in the histograms of Fig. 7. For GLAD and LASER separately, the blended model was able to obtain better single-step velocity prediction when there was more similarity in trained and tested flow features and the training and test drifters were geographically and temporally proximal. Even with the handicap of selectively limiting the kinds of features previously seen by the data-driven component, the blended model was still able to outperform the purely deterministic approach. Analysis of one such trajectory is further discussed in Section 4.

The last dataset that the blended model was tested against was a combination of both LASER and GLAD drifters, referred to as LASER+GLAD in Fig. 7e-f. This increasingly complex collection of drifters spanned both winter and summer conditions in the Gulf with different flow features and disparate evidence of inertial oscillations. Again, training with a random selection of LASER and GLAD drifters allowed the model to outperform the single-step predictive capability of non-repetitive test method, and both blended model approaches outperformed the purely deterministic model. There still remain large RMSE values for the combined blended models, indicating an inability to always make the proper corrections, or recognize sub-mesoscale features from a wider library of options.

For all three datasets, training with a randomly selected training set resulted in the lowest RMSE for single-step predictions. As the LSTM was trained with single-step predictions, and because of the greater similarity in training and test drifter trajectories for random training, this improvement supports the notion that a blended model can increase model accuracy around flow features similar to those seen previously. Because of the currently limited amount of high-resolution training data available globally, a more realistic forecasting situation was simulated by multi-step predictions for non-repetitive training on the LASER training set. Surprisingly, we found similar skill score improvements to the random training, when test and training drifters were often from the same deployment cluster. Figure 8 shows analogous skill score data as Figure 4 for the non-repetitive training. Again, many of the lowest skill scores from the reduced MR predictions were avoided with the blended model, and the blended model outperformed for 63\% of the test set. This is a slight decrease from 67-68\% from the random training and the range of skill scores was larger for the non-repetitive case, actually extending into negative values. Because we divided the training and test drifters dataset differently, designating all of the negative skill score drifters from Figure 5 as training drifters, the strong negative scores from the reduced-MR drifters from Figure 4 are no longer present.

\section{Discussion}

A blended deterministic and neural network model for predicting Lagrangian ocean drifter trajectories substantially outperformed a deterministic reduced-order Maxey-Riley model for two experimental datasets in the Gulf of Mexico. We approximated drifter velocities with significantly less error (Fig. 7) along drifter paths, indicating a better representation of the underlying flow field is available if given the appropriate initial conditions. The blended model also had increased multi-step predictive capabilities, outperforming the ability of the reduced MR equations to correctly predict drifter locations at timescales up to two months. Skill scores were typically higher for both repetitive and non-repetitive training, with improved behavior linked to correct identification of dominant flow features in several cases.

Specifically, the blended model was also able to accurately trace a large vortical feature in Fig. 3 that was missed by the Maxey-Riley equation. Though this feature was mesoscale in size, inaccuracies in the baseline model combined with low-resolution flow data prevented resolving this dominant flow feature for several drifters. For oil spill, search and rescue, and ocean ecosystem research purposes, it is precisely these transport features that are of utmost importance to model and resolve.

Model improvement was also evident in the shallow waters near the outlet of the Mississippi River where a large number of the LASER drifters remained for several weeks, but currents in the geostrophic data brought the MR drifters off of the continental shelf. In this region, there is complex fresh-saltwater mixing as previously highlighted by \citet{Goncalves2019}, and strongly influential sub-mesoscale signatures not immediately evident in the geostrophic currents or 10-meter wind fields. The colocation of poorly performing drifters and the mixing interface was verified in satellite sea surface temperature measurements during the LASER campaign. Near this transport barrier, the MR model obtained negative skill scores whereas the randomly trained blended model performed just as well as in other areas. As the two models only differed by the inclusion of the LSTM component, this improvement can be attributed to the ability of the neural network to mimic a real drifter velocity from a similar signature in  $\boldsymbol{\xi}(t)$ that was seen during training. It is also worth reiterating that $\boldsymbol{\xi}(t)$ contains no position information, and thus the model was able to make changes based solely on the time history of the physical motion of the drifter and the surrounding fluid flow.

Testing the blended model with more difficult tasks, such as using distinctly dissimilar training and test trajectories or testing over a broader range of dynamics (combined summer and winter) started to highlight the current limitations of a blended-model approach. While the blended-model always outperformed the baseline deterministic model upon which it was improving (Fig. 7), the improvements were less significant than with random training. The incremental model improvement from MR equations, to non-repetitive training, to random training is examined in Figure 9. Single step predictions for the drifter represented by the red dot in Fig. 6 is examined for the three models in Fig. 9 (top). Clearly, the modeled velocity approaches the drifter velocity as the data-driven component is introduced and improves. 

This single-step improvement is complemented by trajectories prediction improve in Fig. 9 (bottom row). Because the MR model under-predicts the drifter speed, the ensemble trajectories quickly diverge from the real drifter path, resulting in a poor skill score (inset). The non-repetitive model performs with enhanced physics, resulting in an improved skill score, even though the nearest trained drifter was deployed nearly 70 km and 4 hours away. The random training provided a training drifter deployed only 34 km away, and the arcing features of the real drifter were most closely mimicked. An increased proximity to training drifters likely aided in the improved model performance, but this difference is largely missed when looking at minimal differences in cross correlation and MS coherence in Figure 6.

It is still an open problem to quantitatively describe how Lagrangian trajectory similarity influences this kind of machine learning. As seen in Fig. 6 and 9, neither correlation nor spectral coherence between a tested drifter and the training set can fully account for the differing performance between the random and non-repetitive training methods. Higher correlation and spectral coherence values typically result in lower RMSE values for single-step prediction, but low RMSEs are also possible with poor correlation or coherence. While using random training to obtain similar deployment coordinates (space and time) for test and training drifter results in better blended model performance, this does little to help us understand why drifters seeded away from the training set can also perform well, or to predict if a blended Lagrangian model will perform well in a region without drifters nearby. Further investigation of a definition of Lagrangian similarity for this kind of forecasting and modeling would be beneficial for constraining models to appropriate domains of application, but is currently beyond the scope of this initial investigation as there is no guarantee of stationarity with this kind of data.

\section{Conclusions}

Because of the lack of availability of high spatial- and temporal-resolution flow data that can resolve sub-mesoscale motions, the use of machine learning for improving Lagrangian drifter models shows great promise. The blended model approach pursued here prevents user-defined parameterization and tuning with respect to one specific drifter geometry, drag, buoyancy corrections and other factors that can complicate physical models of drifter trajectories. For drifters of the same type in similar environments and similar seasons, signatures in drifter velocity and acceleration time series can be identified with recurrent neural networks and used to auto-correct velocity fields, even in the presence of strong coastal temperature and salinity gradients, and inertial oscillations. With the current availability of high-resolution drifter data, it is difficult to evaluate the utility of a one-size-fits-all blended model for all seasons. Seasonally restricting the domain of training and testing, however, provides significant improvements over deterministic drifter modeling.

Similar to other machine learning approaches, this blended model may be further improved by training on additional data sets that include more sub-mesoscale features. This added physical complexity will likely need to be met with increased LSTM complexity. The flexible blended approach is also adaptable to include additional physical parameters in its training (e.g. temperature, salinity) should a physical basis for their inclusion be justified for a different baseline model. As well, the blended approach can be adapted for further advances in deterministic drifter modeling, such as the recent slow-manifold developments of Beron-Vera et al. (2019). The broad implications of machine learning improvements could be particularly advantageous to the oceanography community as data from future field campaigns can be harnessed to better understand the motions present in a particular region, and compared to others. Because the blended approach uses highly accurate Lagrangian data to approximate the underlying physics of drifter transport, we believe there is a theoretical finite limit on additional training datasets that are necessary for accurate drifter modeling in all regions of the globe. Applications of this method to other regions are currently underway.








%
\acknowledgments
TPS acknowledges support from a Doherty Career Development Chair and the ONR MURI grant N00014-17-1-2676. All data required to replicate the analysis in this manuscript is freely available online. The LASER experiment ocean drifter trajectory data is hosted by the Gulf of Mexico Research Initiative (DOI: 10.7266/N7W0940J). ECMWF ERA-5 hourly 10m wind fields are available from The Copernicus Programme Climate Change Service (https://climate.copernicus.eu). The AVISO geostrophic current velocity product used in this study, “Global Ocean Gridded L4 Sea Surface Heights and Derived Variables Reprocessed”, is freely available and hosted by the Copernicus Marine Environment Monitoring Service (http://marine.copernicus.eu).

%
\appendix

The Maxey-Riley equation define the path {\bf x}(t) of a small spherical particle in a time-dependent flow as follows:

\begin{equation}
\begin{split}
\dot{{\bf x}}=&{\bf v}({\bf x},t) \\
\rho_p \dot{{\bf v}} = &\rho_f \frac{D{\bf u}}{Dt} + (\rho_p -\rho_f) {\bf g} \\
&-\frac{9 \nu \rho_f}{2a^2} ({\bf v} - {\bf u} - \frac{a^2}{6}\Delta {\bf u}) \\
& - \frac{\rho_f}{2}(\dot{{\bf v}} - \frac{D}{Dt}({\bf u}+\frac{a^2}{10}\Delta{\bf u})) \\
& -\frac{9 \nu \rho_f}{2a} \sqrt{\frac{\nu}{\pi}}\int_0^t\frac{1}{\sqrt{t-s}}(\dot{{\bf v}}(s)-\frac{d}{ds}({\bf u}+\frac{a^2}{6}\delta{\bf u}))ds \\
\end{split}
\end{equation}

where the individual force terms on the right hand side of Eq A1 are the force of the fluid on the particle by the undisturbed flow, the buoyancy force, the Stokes drag, an added mass term from the fluid moving with the particle, and the Basset-Boussinesq memory term. We consider four additional influences on an ocean drifter trajectory following the work the suggestions of \citet{Niiler1987}, \citet{Geyer1989}, \citet{Edwards2006} and  \citet{Beron-Vera2016}. These forces include the wind drag on floats above the sea surface, drag on the tether and drogue induced by shear currents, motions induced by surface waves, and the Coriolis Force.

The wind drag on the exposed float and Coriolis effect can be considered directly, but explicit forms of the shear current drag and surface wave influences are beyond the scope of this research remain open problems in oceanography. Near-surface shear current effects are significantly complicated by stratification in the ocean caused by temperature or salinity gradients, and the degree of turbulence in a specific ocean region \citep{Niiler1995}. It is assumed that the Ekman effect decreases exponentially over the Ekman layer depth, $D_E$ (estimated at 10-40 m for this study) in a well-mixed upper ocean layer. The current changes from $\mathbf{u}=\mathbf{u}_e+\mathbf{u}_g$ at the surface to $\mathbf{u}=\mathbf{u}_g$ at a depth of $D_E$ \citep{Stewart2008}. Thus the resulting the shear acting on the drogue and tether is a function of $\frac{du}{dz}$ which is a function of $\frac{u_e}{zD_E}$, depending on stratification and the shape of the exponential. The effect of this motion is thus included as an unknown function of Ekman current, a linear function of windspeed, that will be learned by the machine learning corrections. For a summary of wind-induced wave effects and specific details of wave effects on drifters from this study, refer to Niiler et al. [1987] and Haza et al. [2018], respectively. The effect of this surface wave motions will be included as an unknown function of Ekman current, a linear function of windspeed, and the surface ocean velocity and will also be learned by the neural network.

We then write an adapted form of Eq. A1 for ocean drifter motion as follows:

\begin{equation}
\begin{split}
\rho_p \dot{{\bf v}} = &\rho_f \frac{D{\bf u}}{Dt} \\
& +(\rho_p -\rho_f) {\bf g} \\
&-\frac{9 \nu_f \rho_f}{2a^2} ({\bf v} - {\bf u} - \frac{a^2}{6}\Delta {\bf u}) \\
&- \frac{\rho_f}{2}(\dot{{\bf v}} - \frac{D}{Dt}({\bf u}+\frac{a^2}{10}\Delta{\bf u})) \\
& -\frac{9 \nu_f \rho_f}{2a} \sqrt{\frac{\nu_f}{\pi}}\int_0^t\frac{1}{\sqrt{t-s}}(\dot{{\bf v}}(s)-\frac{d}{ds}({\bf u}+\frac{a^2}{6}\delta{\bf u}))ds \\
& - \nu_a\rho_a\alpha({\bf v} - {\bf u_{wind}}) - f(\rho_p {\bf v} - \rho_f {\bf u})^{\perp} + F({\bf u_ e},{\bf u}) 
\end{split}
\end{equation}

where  $\alpha$ is a drag coefficient depending on the unique geometry and surface area of the float not submerged in the ocean at any point in time, $f$ is the Coriolis parameter and $F$ is an unknown function combining surface wave effects and shear influences. We can remove the Fauxén terms because $\frac{a}{L} \ll 1$, where $L$ is a dominant length scale of our ocean flow \citep{Wan2018a}.  If we can assume $\frac{a}{\sqrt{\nu_f}}$ is very small, then we can remove the Basset-Boussinesq term. After rescaling space, time, and velocity by a characteristic length scale $L$, characteristic time scale $T=L/U$ and characteristic velocity $U$, we have the 

\begin{equation}
\begin{split}
\dot{{\bf v}} -\frac{3R}{2}\frac{D{\bf u}}{Dt} = -\mu({\bf v} - {\bf u}) +\left(1-\frac{3R}{2}\right){\bf g} -R\left[f(\delta_p{\bf v}-{\bf u})^{\perp}+\delta_a\nu_a\alpha({\bf v}-{\bf u_{wind}})-F({\bf u_e},{\bf u})\right]
\end{split}
\end{equation}

with $R=\frac{2\rho_f}{\rho_f+2\rho_p}$, $\mu=\frac{R}{St}$, $St=\frac{2}{9}\big (\frac{a}{L}\big)^2 Re$, $\delta_a=\frac{\rho_a}{\rho_f} \ll 1$, $\delta_p=\frac{\rho_p}{\rho_f} \approx 1$ by design and $R$ is close to $\frac{2}{3}$. Introducing the small parameter $\epsilon = \frac{1}{\mu} \ll 1$, we can rewrite our system as

\begin{align*}
&{\bf \dot{x}} = {\bf v} \\
&\epsilon\dot{{\bf v}}({\bf x},t)= {\bf u}({\bf x},t) - {\bf v} + \epsilon\left[ \frac{3R}{2}\frac{D{\bf u}({\bf x},t)}{Dt} + \left(1-\frac{3R}{2}\right){\bf g} -R\left[f(\delta_p{\bf v}-{\bf u})^{\perp}+\delta_a\nu_a\alpha({\bf v}-{\bf u_{wind}})-F({\bf u_e},{\bf u})\right]\right] \\
\end{align*}

Here, assuming the form of the Ekman current velocity ${\bf u_e}=\frac{0.0127}{\sqrt{sin|\lambda|}} \mathbf{A} {\bf u_{wind}}$ where $\mathbf{A}$ is a 45 degree rotation to the right in the Northern hemisphere, and $\lambda$ is latitude, we combine the additional slow wind functions into one function of Ekman velocity and ocean current, and drifter velocity $H({\bf u_e}({\bf x},t),{\bf u},{\bf v})$:

\begin{equation}
\begin{split}
\epsilon\dot{{\bf v}}({\bf x},t)= {\bf u}({\bf x},t) - {\bf v} + \epsilon\Big [ \frac{3R}{2}\frac{D{\bf u}({\bf x},t)}{Dt} + (1-\frac{3R}{2}){\bf g} -Rf(\delta_p{\bf v}-{\bf u})^{\perp}+H({\bf u_e}({\bf x},t),{\bf u},{\bf v})\Big ].
\end{split}
\end{equation}

Upon introducing the fast-time by letting $\epsilon \tau = t - t_0$, $\phi = t_0 + \epsilon \tau$ be a dummy variable, and denoting differentiation with respect to $\tau$ by prime, we can rewrite our system as the autonomous dynamical system

\begin{equation}
\begin{split}
&{\bf x'} = \epsilon{\bf v} \\
&\phi' = \epsilon \\
& {\bf v}'({\bf x},\varphi)= {\bf u}({\bf x},\phi) - {\bf v} + \epsilon\Big [ \frac{3R}{2}\frac{D{\bf u}({\bf x},\phi)}{Dt} + (1-\frac{3R}{2}){\bf g} -Rf(\delta_p{\bf v}-{\bf u})^{\perp}+H({\bf u_e}({\bf x},\phi),{\bf u},{\bf v})\Big ]
\end{split}
\end{equation}

The $\epsilon=0$ limit of the system

\begin{equation}
\begin{split}
&{\bf x'} = 0 \\
&\phi' = 0 \\
& {\bf v}'({\bf x},\varphi)= {\bf u}({\bf x},\phi) - {\bf v}.
\end{split}
\end{equation}

has a set of fixed points. Taylor expanding these solutions in $\epsilon$ gives

\begin{equation}
M_\epsilon=\left\{ ({\bf x},\phi,{\bf v})  \; {\bf  : } \;  {\bf v}=  {\bf u}({\bf x},\phi) +\epsilon{\bf u}^1({\bf x},\phi) +\cdots + \epsilon^r{\bf u}^r({\bf x},\phi) + O(\epsilon^{r+1}), ({\bf x},\phi)\in D_0   \right\}.
\end{equation}

Equation (6) restricted to $M_\epsilon$ is a slowly varying system of the form 

\begin{equation}
\begin{split}
{\bf x}' &=  \epsilon{\bf v}|_{M_\epsilon} \\
& = \epsilon[{\bf u}({\bf x},\phi) +\epsilon{\bf u}^1({\bf x},\phi) +\cdots + \epsilon^r{\bf u}^r({\bf x},\phi) + O(\epsilon^{r+1})]
\end{split}
\end{equation}

Now, differentiating 

\begin{equation}
{\bf v} = {\bf u}({\bf x},\phi) +\sum_{k=1}^r \epsilon^k{\bf u}^k({\bf x},\phi) + O(\epsilon^{r+1})
\end{equation}

with respect to $\tau$ gives 

\begin{equation}
{\bf v}' = {\bf u_x x}'+{\bf u}_\phi\phi' +\sum_{k=1}^r \epsilon^k\left[{\bf u}^k_{\bf x}{\bf x}'+{\bf u}^k_\phi\phi'  \right] + O(\epsilon^{r+1}).
\end{equation}

From Equations (6) and (10) we also have

\begin{equation}
\begin{split}
{\bf v}' = -\sum_{k=1}^r \epsilon^k{\bf u}^k({\bf x},\phi)+\epsilon\Big [ \frac{3R}{2}\frac{D{\bf u}({\bf x},\phi)}{Dt} &+ (1-\frac{3R}{2}){\bf g} -Rf(\delta_p\sum_{k=0}^r \epsilon^k{\bf u}^k({\bf x},\phi)-{\bf u})^{\perp}\\
&+H\left({\bf u_e}({\bf x},\phi),{\bf u},\sum_{k=0}^r \epsilon^k{\bf u}^k({\bf x},\phi)\right)\Big ]+O(\epsilon^{r+1}).
\end{split}
\end{equation}

Equating terms of equal power of epsilon in (11) and (12) gives a leading order approximation in our original time, $t$

\begin{equation}
\dot{{\bf x}} = {\bf u}({\bf x},t)+\epsilon\Big [ \left(1-\frac{3R}{2}\right)\left(\frac{D{\bf u}({\bf x},t)}{Dt} - {\bf g}\right) -Rf(\delta_p-1){\bf u}^{\perp}+H({\bf u_e}({\bf x},t),{\bf u})\Big ]
\end{equation}

The sign of the Coriolis term is the same as in \citet{Beron-Vera2015}, but appears different because of multiplication by $R$ and differing definitions $\delta$ and $\delta_p$. Without better models of the wind-wave coupling, upper-ocean stratification, or more information on drifter specifics (size, geometry, density, compressibility, float-drifter tension, etc.), we have limited our modelling approach to a leading order approximation of drifter behavior, and let a machine learning algorithm learn the remaining terms as introduced by \citet{Wan2018} and \citet{Wan2018a}. We remove gravity and are left with the following form our blended model:

\begin{equation}
\begin{split}
\dot{{\bf x}} &=  {\bf v}({\bf x},t) + \epsilon\left[\left(\frac{3R}{2} - 1 \right)\frac{D{\bf u}}{Dt} -Rf(\delta_p-1){\bf u}^{\perp}\right]+\mathbf{G}(\boldsymbol{\xi}(t),\boldsymbol{\xi}(t-\tau),\boldsymbol{\xi}(t-2\tau), \dots) \\
 \boldsymbol{\xi}(t) &=  \left[{\bf u}({\bf x},t),{\bf u_e}({\bf x},t),  \frac{D{\bf u}}{Dt} \right]
 \end{split}
\end{equation}

where $\tau$ is a time delay. We are thus forcing the data-driven model $\mathbf{G}$ to learn our as yet unknown float-drifter-wind coupling function $H$, as well as the higher order terms from the Taylor expansion.






%
%
%
 \bibliographystyle{ametsoc2014}
\bibliography{references}

%
\begin{table}[t]
\caption{Mean values of velocity statistics (Correlation and Mean MS coherence) comparing test and training datasets for LASER and GLAD drifter trajectories. Values outside parentheses are for randomly selected training sets, and in parentheses are for ordered, non-repetitive dataset}
\begin{center}
\begin{tabular}{ccccrrcrc}
\hline\hline
$ $ & $Laser (Winter)$ & $GLAD (Summer)$ \\
\hline
 Max Correlation & 0.8 (0.7) & 0.8 (0.74)  \\
 Max Mean MS Coherence & 0.56 (0.56) & 0.8 (0.62) \\
 Dist. to Trained Drifter (km) & 5.2 (19.1) & 4.6 (13.7)  \\
  Time to Trained Drifter (hrs) & 0.95 (18.72) & 0.25 (79) \\
\hline%
\end{tabular}
\end{center}
\end{table}

%

\begin{figure}[t]
  \noindent\includegraphics[width=39pc,angle=0]{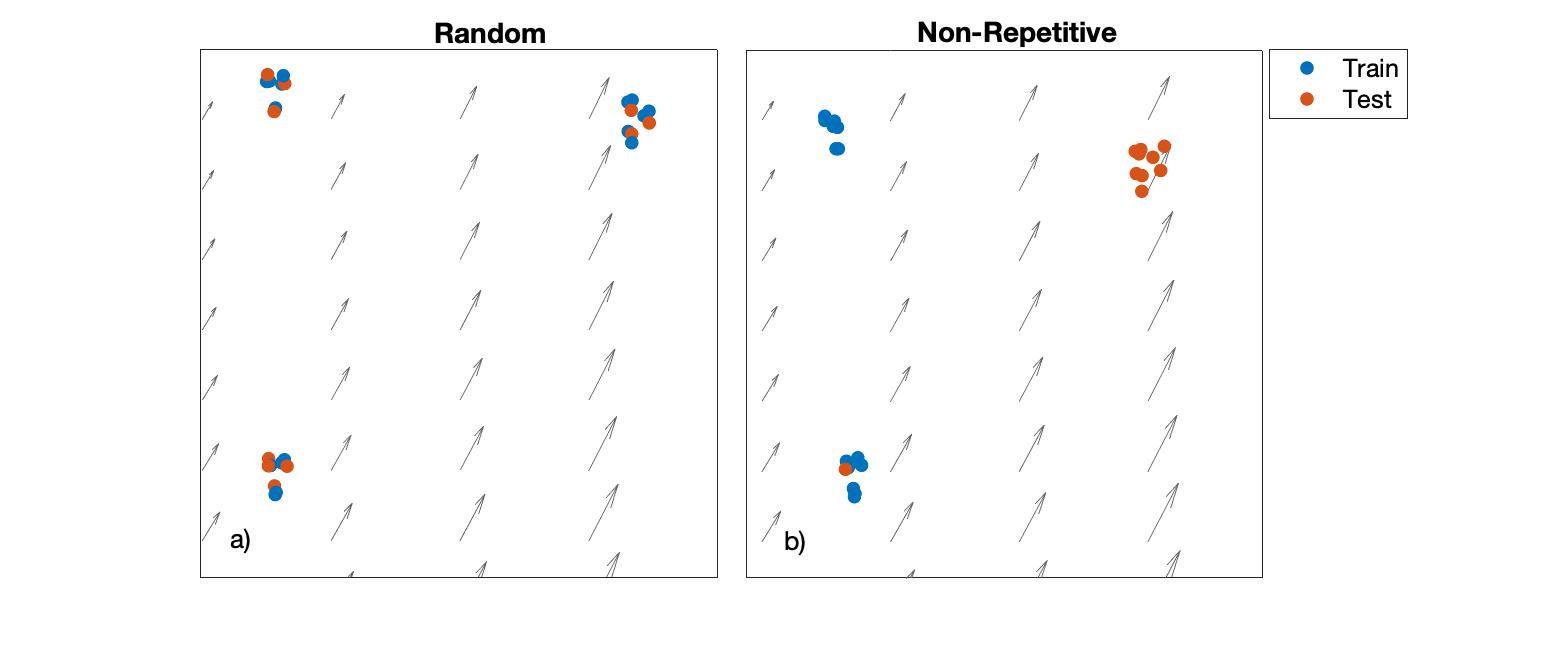}\\
  \caption{Two types of training and test drifter selection used in this study. Random selection maximizes training and test drifter trajectory similarity while non-repetitive training provides a more difficult test for blended-model performance by limiting training and test set similarities. Drifter initial positions are overlain on ocean current vectors.}
\end{figure}

\begin{figure}[t]
  \noindent\includegraphics[width=39pc,angle=0]{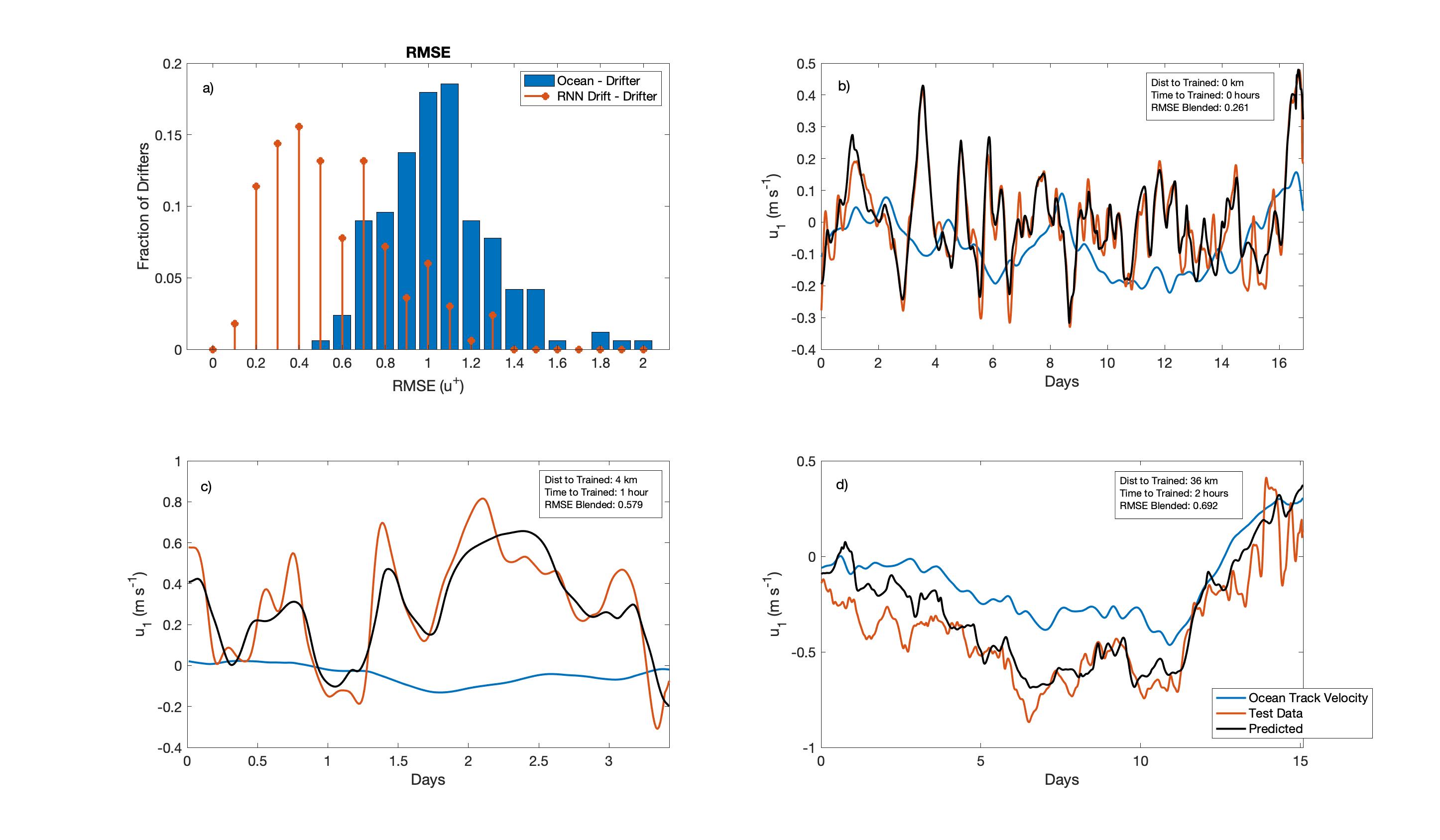}\\
  \caption{Single-step prediction velocity test comparisons. Subplot (a) shows histograms of RMSE for non-dimensionalized zonal components $(u^+=\frac{u_1}{\bar{\mathbf{u}}})$ for the reduced Maxey-Riley equations and the blended model. Subplots (b-d) shows examples of real drifter velocities, blended model velocities and the underlying low-resolution ocean current velocity. A significant reduction in error can be found with the blended approach.}
\end{figure}

\begin{figure}[t]
  \noindent\includegraphics[width=39pc,angle=0]{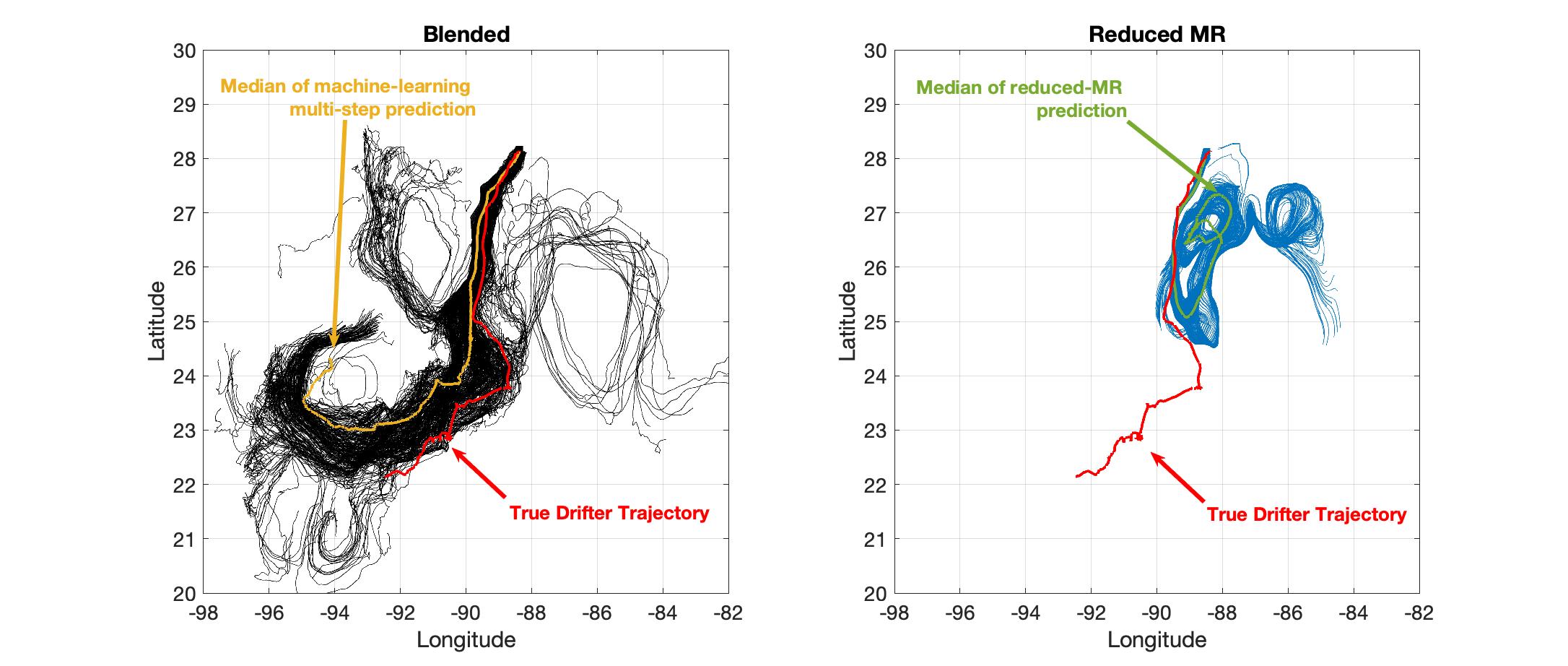}\\
  \caption{Multi-step ensemble trajectory prediction comparisons. The left subplot shows an ensemble of blended model drifters in black, and their median trajectory in yellow. The right subplot shows an ensemble of reduced-MR drifters and their median trajectory in green. On both plots, the true drifter trajectory is indicated by the bold red line.}
\end{figure}

\begin{figure}[t]
  \noindent\includegraphics[width=39pc,angle=0]{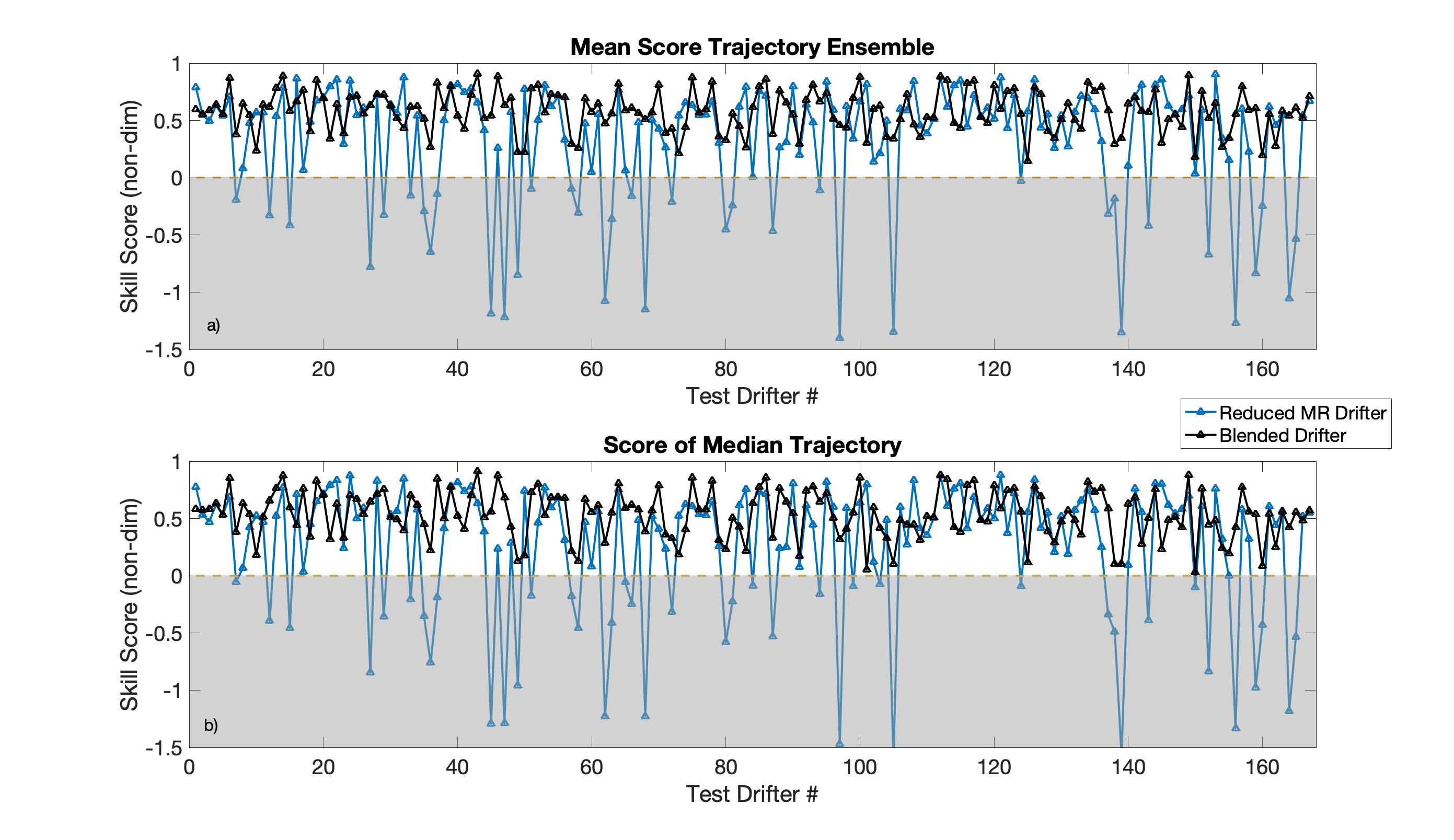}\\
  \caption{Skill scores of the randomly trained blended model and MR model on the LASER dataset. Subplot (a) shows the mean skill score of all drifters in an ensemble prediction for each LASER test drifter. Subplot (b) shows skill scores for the median trajectory from an ensemble prediction for LASER drifters. The gray shaded area indicates the region where model and real drifter trajectory separation outpaces the actual drifter motion.}
\end{figure}

\begin{figure}[t]
  \noindent\includegraphics[width=39pc,angle=0]{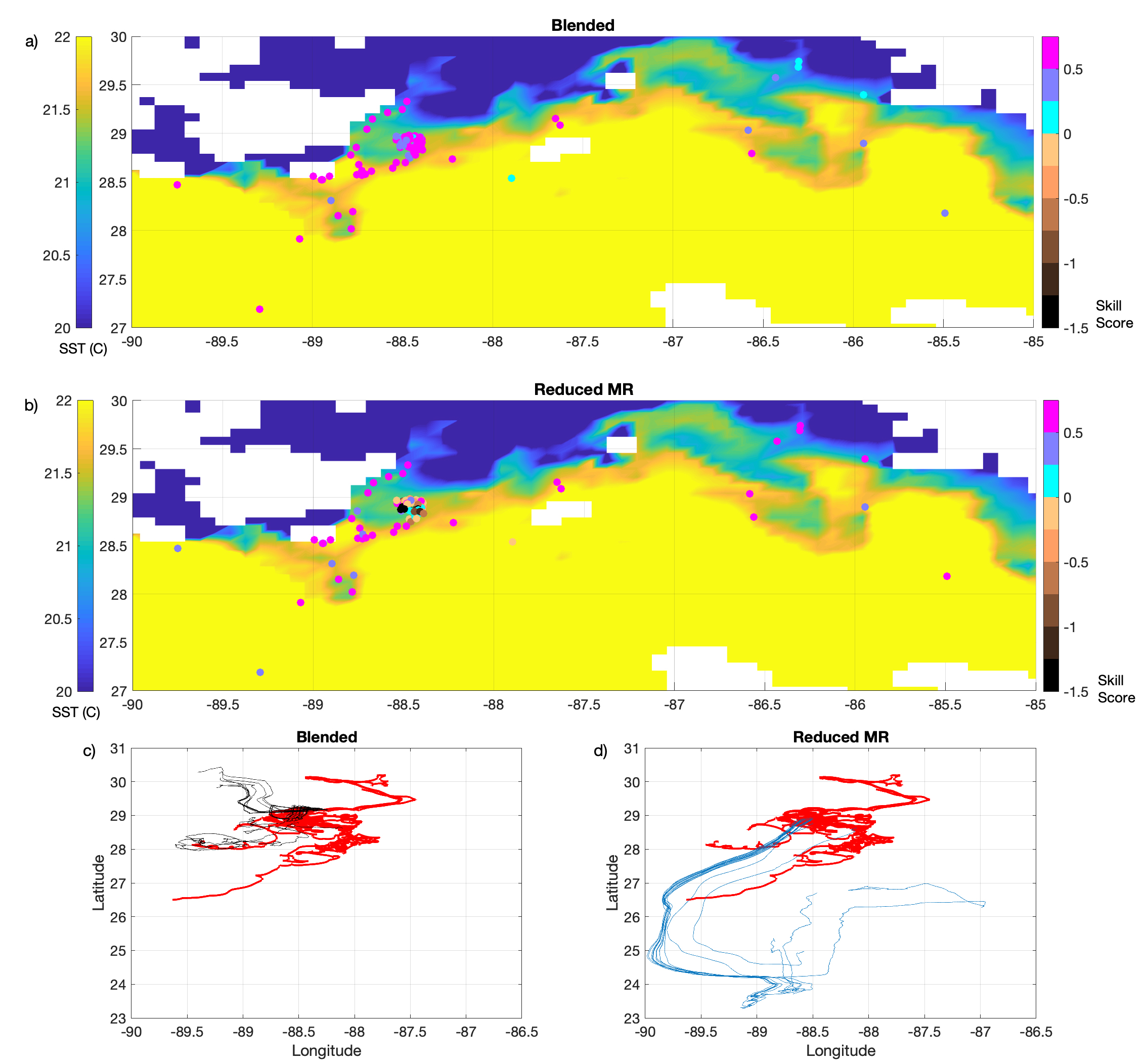}\\
  \caption{Positions and skills scores for the blended and reduced MR models on February 7, 2016 overlaid on daily Aqua MODIS Thermal IR sea surface temperature measurements. White areas in the top-left correspond with the outlet of the Mississippi river. Subplots c-d show the median trajectories of the drifters with negative skill scores. Reduced MR drifters quickly leave the coastal zone while the blended model drifters wander for their entire lifetime. The time elapsed for trajectories in c-d is approximately 23 days.}
\end{figure}

\begin{figure}[t]
  \noindent\includegraphics[width=39pc,angle=0]{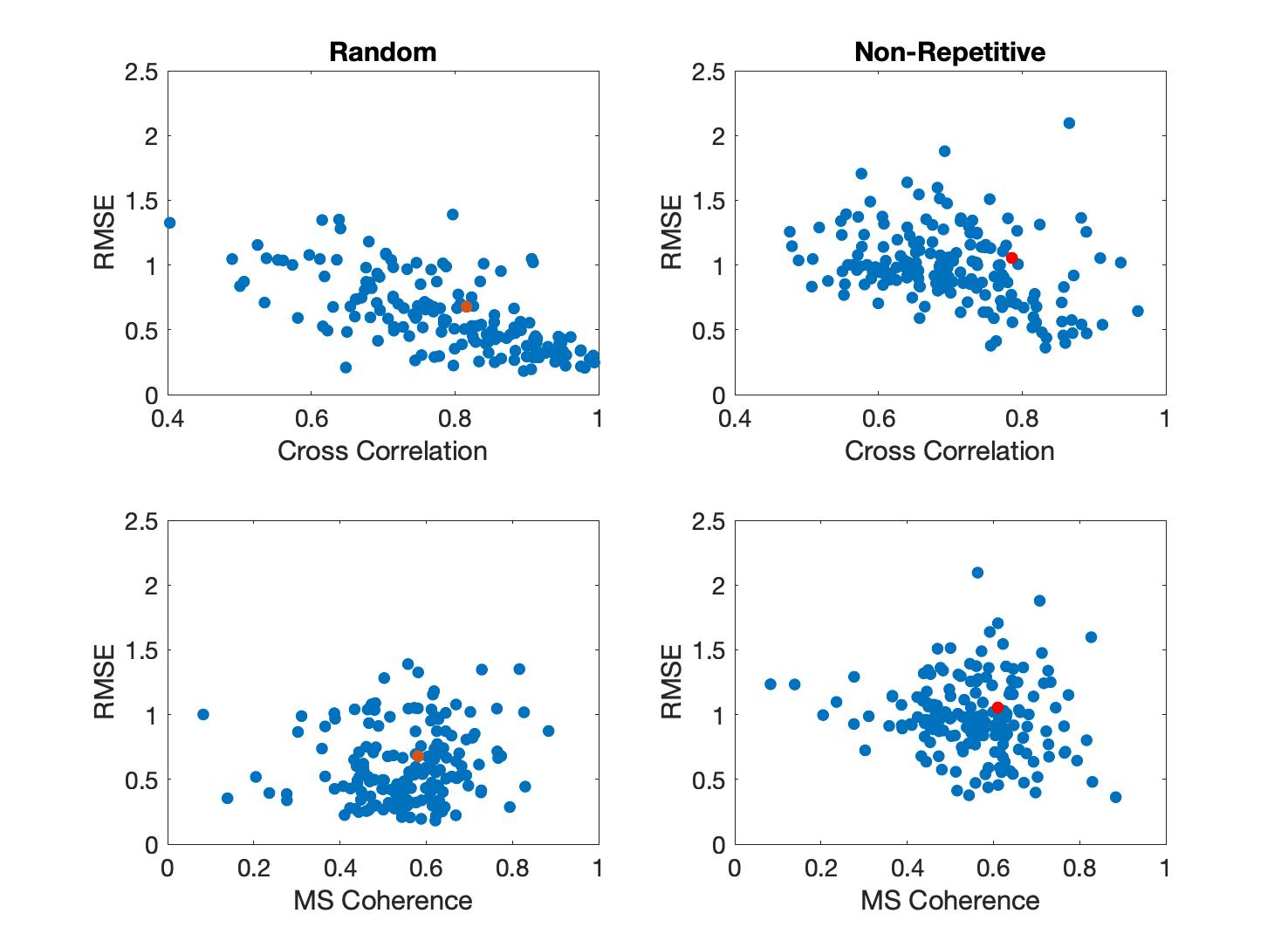}\\
  \caption{Cross Correlation and mean MS coherence vs RMSE values for zonal velocity of blended drifter modeling from random and non-repetitive training methods. There is clearest decrease in RMSE with increased correlation for the random data, while all other relationships remain less defined. The same test drifter is marked by a red dot in both subplots and further explored in section 4.}
\end{figure}

\begin{figure}[t]
  \noindent\includegraphics[width=39pc,angle=0]{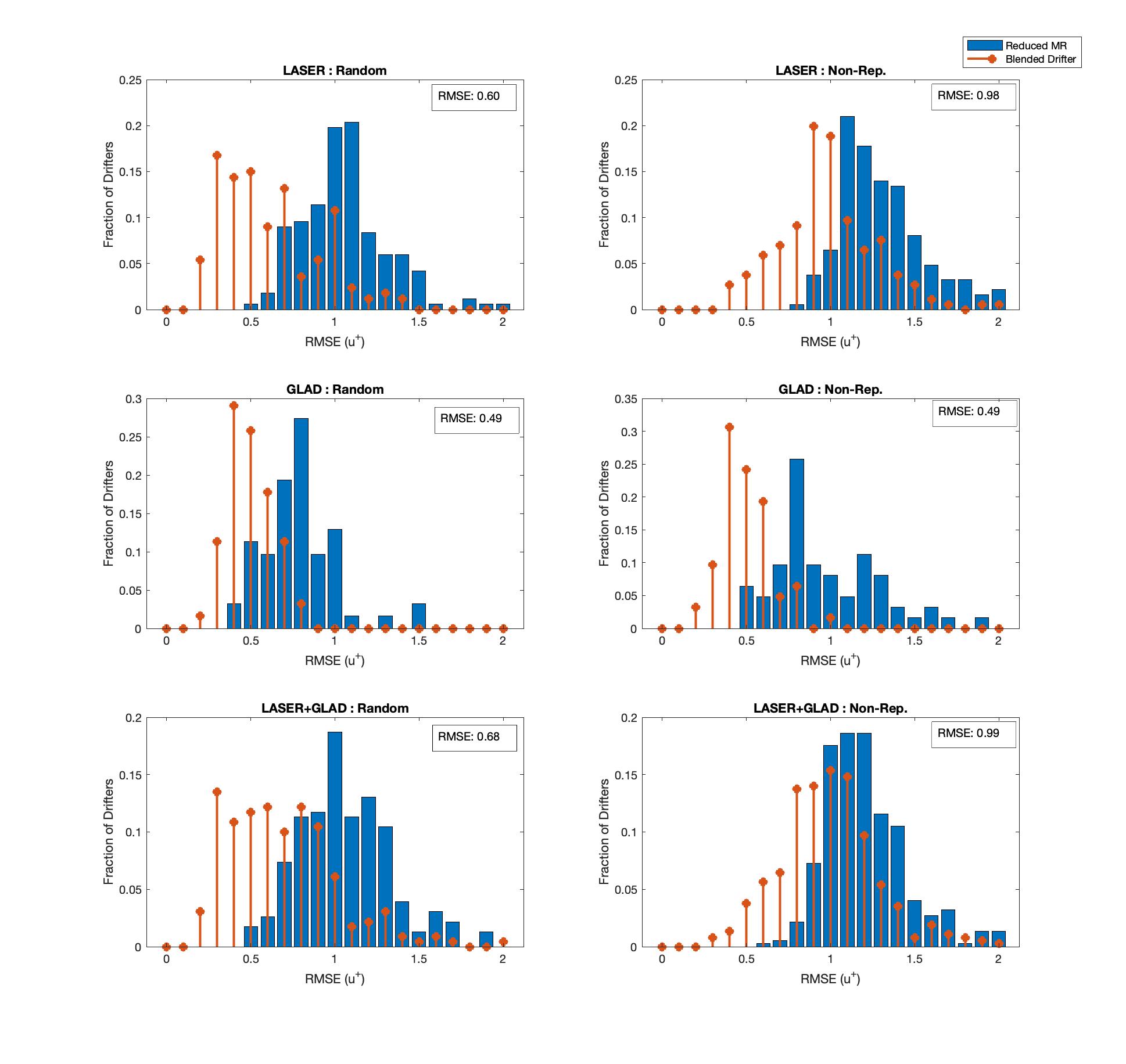}\\
  \caption{RMSE histograms of non-dimensionalized zonal velocity for surface drifters during the LASER (Winter), GLAD (Summer), and LASER+GLAD datasets (rows top to bottom, resp.) for random and non-repetitive trajectory training (left and right, resp.). Blue bars indicate the performance of the deterministic reduced Maxey-Riley equations and orange stem plots represent the performance of the blended model.}
\end{figure}

\begin{figure}[t]
  \noindent\includegraphics[width=39pc,angle=0]{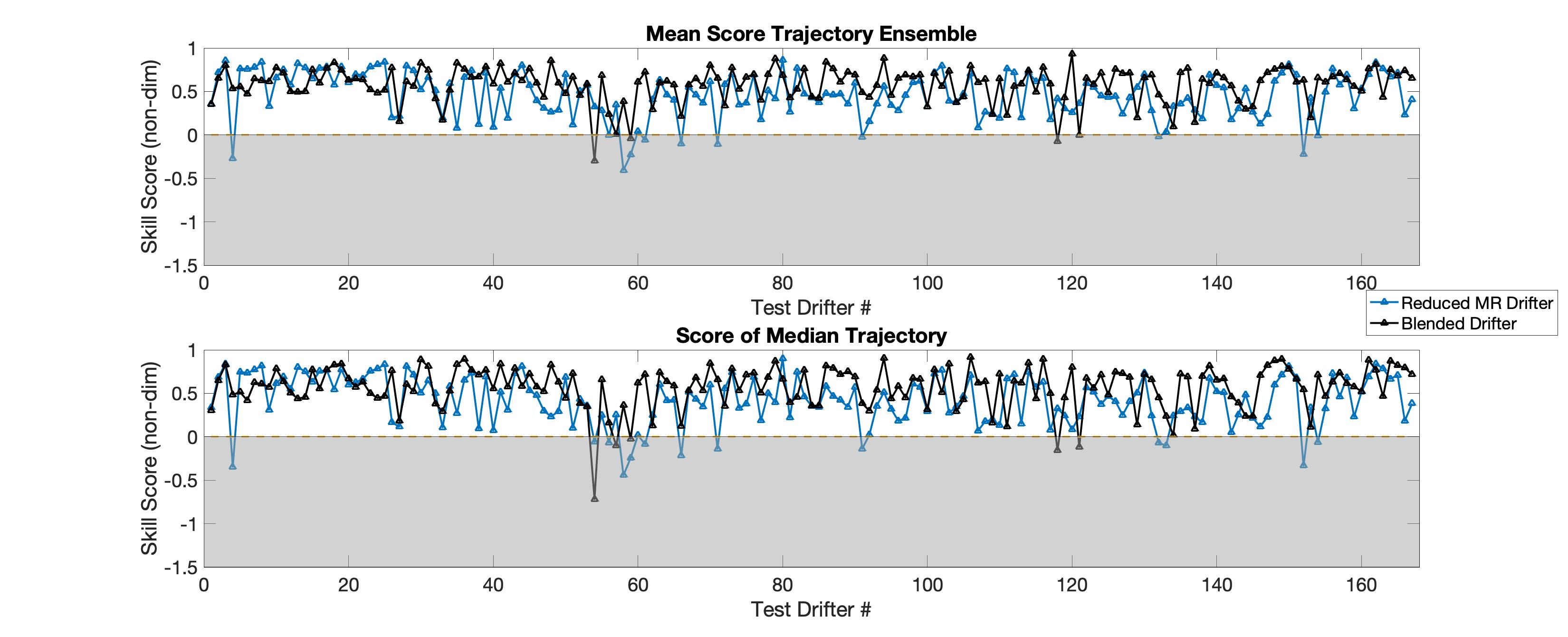}\\
  \caption{Skill scores calculated for the non-repetitive winter (LASER) test and training dataset, similar to Figure 3. For this case, the blended model still outperforms the reduced MR predictions though the difference is less dramatic as the negative skill-score mixing-region drifters highlighted in Fig. 5 are now in the training set.}
\end{figure}

\begin{figure}[t]
  \noindent\includegraphics[width=39pc,angle=0]{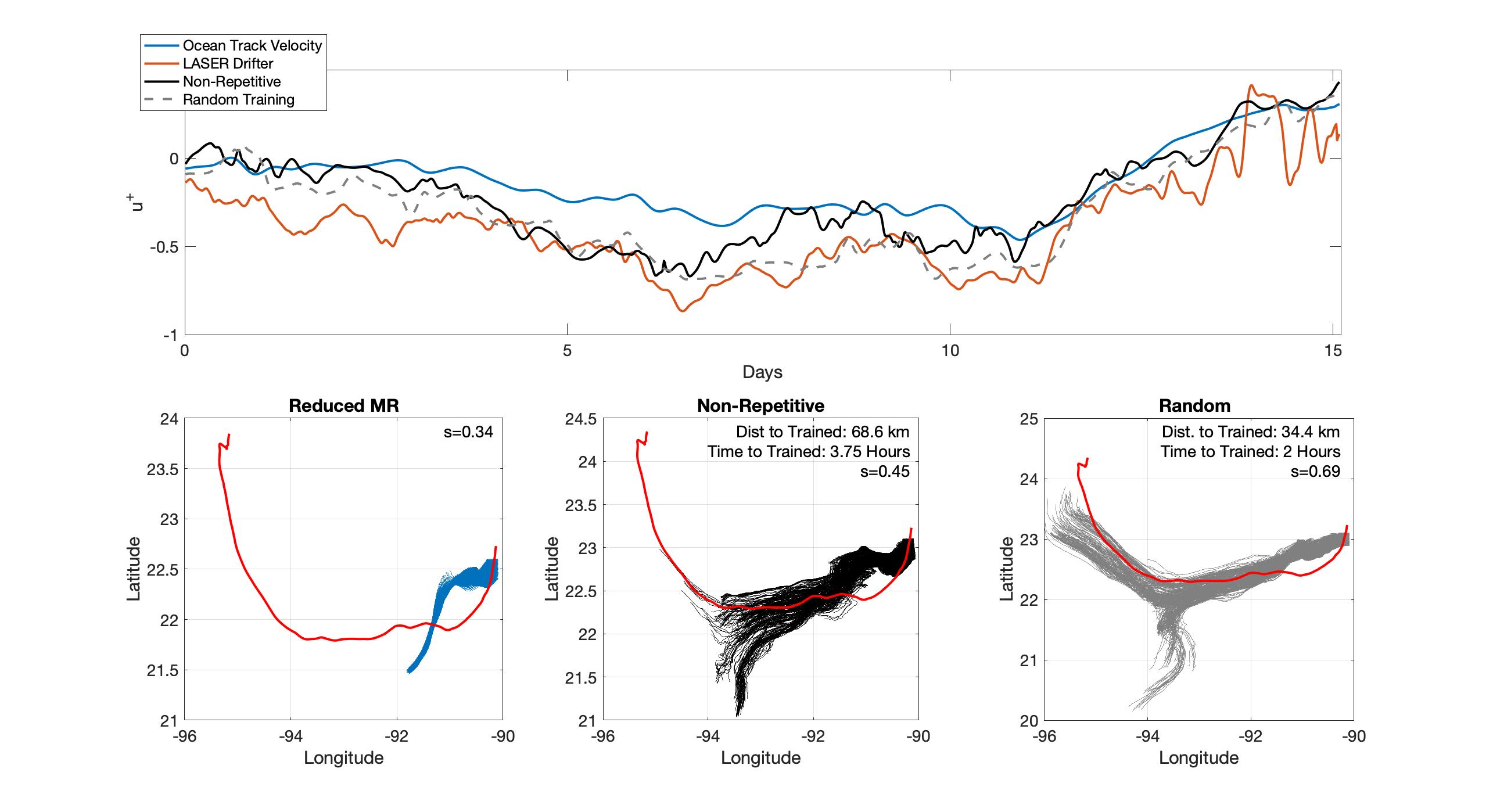}\\
  \caption{Trajectory and velocity time series predictions for the non-repetitive and random training sets, and the reduced MR model, as well as the true drifter trajectory and zonal velocity.}
\end{figure}

\end{document}